\documentclass[11pt,a4paper]{article}
\usepackage[latin2]{inputenc}
\usepackage[english]{babel}
\usepackage{a4wide}
\usepackage[english]{babel}
\usepackage[nottoc]{tocbibind}
\usepackage{amsmath,amsfonts,amstext,amsbsy,amsopn,amssymb,amsthm,amscd}
\usepackage{amsxtra,latexsym}
\usepackage{bbm}
\usepackage{exscale}
\usepackage{paralist}
\usepackage[article]{math_env} 
\numberwithin{equation}{section}
\usepackage{dmath_abbr_wias}
\DeclareMathOperator{\curl}{curl}
\DeclareMathOperator{\sym}{sym}
\usepackage{graphicx,mathrsfs}
\usepackage{psfrag}
\DeclareMathOperator{\Curl}{Curl}
\DeclareMathOperator{\dyw}{div}
\DeclareMathOperator{\Dyw}{Div}
\DeclareMathOperator{\Sym}{Sym}

\newcommand{\nn}{\nonumber}

\renewcommand{\D}{{\mathbb C_{{\rm e}}}}

\newcommand{\Ha}{{\mathbb C_{{\rm micro}}}}
\newcommand{\Lc}{{\mathbb L}_{{\rm c}}}

\newcommand{\di}{{\mathrm d}}

\newcommand{\idd}{\mathbbm{1}}

\usepackage{color}

\definecolor{dred}{rgb}{.8,0,0}
\definecolor{ddmagenta}{rgb}{0.7,0,0.9}
\definecolor{ddcyan}{rgb}{0,0.2,1.0}
\definecolor{Green}{rgb}{0.,0.5,0.}

\newcommand{\DDDS}{\color{blue} }
\newcommand{\DDDE}{\color{black}}

\makeatletter
\let\@fnsymbol\@arabic
\makeatother

\begin{document}
\title{Global regularity for a physically nonlinear version of the relaxed micromorphic model on Lipschitz domains}
\author{Dorothee Knees\thanks{
Institute of Mathematics, University of Kassel, Heinrich-Plett Str. 40, 34132 Kassel, Germany, dknees@mathematik.uni-kassel.de
}\,,
\and
Sebastian Owczarek\thanks{Faculty of Mathematics and Information Science, Warsaw University of Technology, ul. Koszykowa 75, 00-662 Warsaw, Poland, sebastian.owczarek@pw.edu.pl}\,,
\and
Patrizio Neff\thanks{Lehrstuhl f\"ur Nichtlineare Analysis und Modellierung, Fakult\"at f\"ur Mathematik, Universit\"at Duisburg-Essen, Campus Essen, Thea-Leymann Str. 9, 45127 Essen, Germany, patrizio.neff@uni-due.de}
}

\date{ }

\maketitle

\begin{abstract}
In this paper, we investigate the global higher regularity properties of weak solutions for a linear elliptic system coupled with a nonlinear Maxwell-type system defined on Lipschitz domains. The regularity result is established using a modified finite difference approach. These adjusted finite differences involve inner variations in conjunction with a Piola-type transformation to preserve the curl-structure within the matrix Maxwell system. The proposed method is further applied to the linear relaxed micromorphic model.

As a result, for a physically nonlinear version of the relaxed micromorphic model, we demonstrate that for arbitrary $\epsilon > 0$, the displacement vector $u$ belongs to $H^{\frac{3}{2}-\epsilon}(\Omega)$, and the microdistortion tensor $P$ belongs to $H^{\frac{1}{2}-\epsilon}(\Omega)$ while $\Curl P$ belongs to $H^{\frac{1}{2}-\epsilon}(\Omega)$.
\end{abstract}

%



%
\section{Introduction}
\label{s:introduction}
The relaxed micromorphic model is an innovative generalized continuum model designed to characterize size effects and the band-gap behavior of microstructured solids. This model provides effective equations that do not require a detailed consideration of the microstructure. Operating within the micromorphic framework, it combines the classical displacement $u:\Omega\subset\R^3\to\R^3$ with a non-symmetric tensor field  $P:\Omega\subset\R^3\to\R^{3\times 3}$, known as the microdistortion, through the variational problem
 \begin{equation}
\label{Variational}
\begin{split}
\int_{\Omega}&\Big(\big\langle \D\sym(\mathrm{D} u-P),\sym(\mathrm{D} u-P)\big\rangle
+\big\langle\Ha\sym P,\sym P\big\rangle
+ \big\langle\Lc\Curl P,\Curl P\big\rangle\\[1ex]
&\hspace{1ex}-\langle f,u\rangle-\langle M,P\rangle\Big)\,\di x\quad\longrightarrow \quad {\rm min\,\,\,\, w.r.t.\,\,\, }(u,P)\,,
\end{split}
\end{equation}
subject to suitable boundary conditions. This coupling between $u$ and $P$ enables a comprehensive description of the complex behavior of microstructured materials without the need to delve into the specifics of their microstructural details \cite{exprolingalberdi,Neff_unfol_22,Ghibadyn,MADEO2016,NGperspective,Owczghibaneffexist}. The tensor $\Lc$  introduces a size-dependence into the model, meaning that smaller samples exhibit relatively stiffer responses. The existence and uniqueness in the static case can be deduced from the incompatible Korn's inequality \cite{Lewintan_Korn2020,lewintan_neff_2021,NeffLP,neff2012poincare}. The constitutive tensors $\D$, $\Ha$ and $\Lc$ can be determined using novel homogenisation strategies \cite{Neffhomog,homo_2,homo_1,schroder2021lagrange}. As $\Ha\rightarrow +\infty$ approaches infinity, the model's response converges to the linear Cosserat (micropolar) model \cite{Ghiba2022CosseratME}. There are already various analytical solutions available for the relaxed micromorphic model in engineering, as mentioned in \cite{rizzi2021torsion}. The solution can be naturally obtained with the conditions that $u$ belongs to the Sobolev space $H^1(\Omega)$ and $P$ belongs to the Sobolev space $H(\text{Curl};\Omega)$. Therefore, the microdistortion $P$ may exhibit discontinuities in the normal direction. In the context of finite element implementation, standard element formulations can be used for the displacement $u$, while N\'ed\'elec spaces are employed for $P$ to achieve optimal convergence rates \cite{Rizzi_boundary_2022,schroder2021lagrange,SKY2022_primal,sky2021hybrid,Sky_H(symcurl)}.

However, there are cases where it is preferable to bypass the N\'ed\'elec framework and instead work with $H^1(\Omega)$ for the microdistortion tensor $P$. In these situations, it is necessary to clarify beforehand whether the regularity of $P$ allows for an accurate outcome in the numerical approximation. In line with this approach, we will now further investigate regularity in the static case and obtain a higher regularity result for weak solutions of the relaxed micromorphic model. We will focous on global regularity for domains with only Lipschitz boundary, as opposed to only considering interior regularity in \cite{KON23} and global regularity for smooth boundaries \cite{KON2023global}. The problem is formulated on a bounded domain $\Omega \subset \mathbb{R}^3$, and the Euler-Lagrange equations for \eqref{Variational} can be expressed as follows (\cite{exprolingalberdi,Ghibadyn,NGperspective,Owczghibaneffexist,OGNdynamicreg}): given positive definite and symmetric material-dependent coefficient tensors $\mathbb{C}_e: \Omega \rightarrow \text{Lin}(\text{Sym}(3), \text{Sym}(3))$, $\mathbb{C}_\text{micro}: \Omega \rightarrow \text{Lin}(\text{Sym}(3), \text{Sym}(3))$, and $\mathbb{L}_c: \Omega \rightarrow \text{Lin}(\mathbb{R}^{3 \times 3}, \mathbb{R}^{3 \times 3})$, we seek to determine a displacement field $u: \Omega \rightarrow \mathbb{R}^3$ and a non-symmetric microdistortion tensor $P: \Omega \rightarrow \mathbb{R}^{3 \times 3}$ that satisfy the system
\begin{align}
\label{Main}
 0&=\Dyw\Big(\bbC_e\sym(\mathrm{D} u - P)\Big) + f\quad \text{in }\Omega,\nn\\
 0&= -\Curl\Big(\bbL_c\Curl P\Big) + \bbC_e\sym(\mathrm{D} u - P) - \bbC_\text{micro}\sym P + M\quad \text{in }\Omega\,
\end{align}
along with appropriate boundary conditions.

In this article, we aim to expand on the findings of global regularity for the problem described in system \eqref{Main} in the case of a smooth boundary, as referenced in the work \cite{KON2023global}. Our objective is to extend this result to the case of a Lipschitz boundary. To achieve this, we will present our methodology in the context of a physically nonlinear version of the relaxed micromorphic model.

\subsection{A nonlinear version of the relaxed micromorphic model}
The relaxed micromorphic model, as introduced above, is a geometrically linear model: the coupling of displacements $u$ and microdistortion $P$ is additive, and the energy terms are quadratic. Notably, the dislocation density tensor $\Curl P$ enters quadratically. Physical considerations (e.g., for materials with distributed dislocations) suggest that a more realistic framework consists of taking only  $\norm{\Curl P}^1$ in the energy. However, in this case, even the existence of minimizers is an open problem. Therefore, in this contribution, we consider especially the contribution $\alpha\norm{\Curl P}^q$, $1<q<2$, and add a small, quadratic term controlling $\norm{\Curl P}^2$. Then the new Euler-Lagrange equations read:
\begin{align}
\label{Main_nonlinear}
 0&=\Dyw\Big(\bbC_e\sym(\mathrm{D} u - P)\Big) + f\quad \text{in }\Omega,\nn\\
 0&= -\Curl\Big(\bbL_c\Curl P+q\,\alpha\norm{\Curl P}^{q-2}\Curl P\Big) + \bbC_e\sym(\mathrm{D} u - P)\\
 &\hspace{4ex}- \bbC_\text{micro}\sym P + M\quad \text{in }\Omega\,.\nn
\end{align}
For such a nonlinear problem, the standard methods of functional analysis \cite{KON2023global,Pauly_Sch_P_1,Pauly_Sch_P_2,KuhnPauly} do not apply any more and we need to resort to the method of inner variations, already employed in \cite{KON23,KS12,Neffkneesreg,Nes06}. In fact, this method is very versatile and allows us to treat even more general, nonlinear situations than \eqref{Main_nonlinear}. Therefore, we exploit the generality of the method in considering the geometrically linear, but physically nonlinear model problem
 \begin{equation}
\label{Variational_nonlinear}
\begin{split}
\int_{\Omega} W(\mathrm{D}u,P,\Curl P)-\langle f,u\rangle-\langle M,P\rangle\,\di x\quad\longrightarrow \quad {\rm min\,\,\,\, w.r.t.\,\,\, }(u,P)
\end{split}
\end{equation}
with associated Euler-Lagrange equations
\begin{align}
\label{Main_nonlinear_1}
 0&=\Dyw\mathrm{D}_{\mathrm{D}u}W(\mathrm{D}u,P,\Curl P) + f\quad \text{in }\Omega,\nn\\
 0&= -\Curl\mathrm{D}_{\Curl P}W(\mathrm{D}u,P,\Curl P) -\mathrm{D}_{P}W(\mathrm{D}u,P,\Curl P) + M\quad \text{in }\Omega\,,
\end{align}
and assuming suitable structure conditions on the Lagrangian $W$, see (W1), (W2) and (W3) in Section 3.

The paper is organized now as follows. Initially, we will introduce the concepts utilized in this article, establishing connections with the set $\Omega$ and the relevant function spaces. Additionally, we will recall properties related to the Helmholtz decomposition, particularly in the case of a Lipschitz boundary. Subsequently, we will formulate assumptions pertaining to the energy $W$ and provide two illustrative examples that satisfy these assumptions. In the ensuing section, we will present the main result concerning the regularity of solutions in the examined nonlinear case. During the proof of the main result, a family of admissible inner variations will be proposed. The defined properties of this family will enable the identification of a suitable test function for the equation involving the microdistortion tensor $P$. A modified Piola transformation for $P$ will be employed, leading to the conclusion that the solutions to problem \eqref{Main_nonlinear_1} exhibit the following regularity: for arbitrary $\epsilon>0$  
$$ u\in H^{\frac{3}{2}-\epsilon}(\Omega)\,,\qquad\, P\in H^{\frac{1}{2}-\epsilon}(\Omega) \,,\qquad\,  \Curl P\in H^{\frac{1}{2}-\epsilon}(\Omega)\,,$$
where the space $H^s(\Omega)$ is understood as the Sobolev-Slobodeckij space of order $s$. For readers interested in alternative methods to demonstrate the regularity of solutions to nonlinear Maxwell type equations, we recommend, for instance \cite{AMMARI200351,Pan_nonlinear_maxwell,SPITZ20195012}.

\section{Background from function space theory}
Before we formulate the nonlinear problem associated with system \eqref{Main}, we will need definitions, assumptions and some general knowledge related to suitable function spaces.

For vectors $a,\,b\in\R^3$, we define the scalar product $\langle a,b\rangle:=\sum_{i=1}^3 a_ib_i$, the Euclidean norm $\abs{a}^2:=\langle a,a\rangle$ and the dyadic product $a\otimes  b=(a_ib_j)_{i,j=1}^3\in\R^{3\times 3}$, where $\R^{3\times 3}$ will denote the set of real $3\times 3$ matrices. 
For  matrices $P,Q\in \R^{3\times 3}$, we 
define the standard Euclidean scalar product $\langle P,Q\rangle:=\sum_{i=1}^3\sum_{j=1}^3 P_{ij}Q_{ij}$ and the Frobenius-norm $\|P\|^2:=\langle P,P\rangle$. $P^T\in \R^{3\times 3}$ denotes the transposition of the matrix $P\in \R^{3\times 3}$ and for $P\in \R^{3\times 3}$, the symmetric part of $P$
will be denoted by $\sym P=\frac{1}{2}(P+P^T)\in \Sym(3)$. 

Now we record some facts related to the geometry of the domain $\Omega$. Let $\Omega\subset\R^3$ be a bounded domain. We assume in this paper that the boundary $\partial\Omega$ is Lipschitz continuous, meaning that it can locally be described as the graph of a Lipschitz continuous function, see e.g. \cite{Grisvard} for a precise definition. 
\begin{definition}(Uniform cone property)
\begin{itemize}
 \item[(a)] A set $\calC\subset \R^3$ is an open (infinite) cone if there exists an open and non-empty subset $C\subset\partial B_1(0)$, where $B_1(0)=\{x\in \R^3\, |\,\, \abs{x}_2^2<1\}$ is the unit ball with respect to the Euclidean norm, such that 
 \[
  \calC=\Big\{x\in \R^3\backslash\{0\}\, |\,\, \tfrac{x}{\abs{x}_2} \in C\Big\}\,.
 \]
 Let $\calC\subset\R^3$ be an open cone and $\rho>0$. Then we define $\calC_\rho:=\calC\cap B_\rho(0)$.
\item[(b)] A bounded domain $\Omega\subset \R^3$ has the uniform (exterior or Poincar\'e cone property) cone property if 
there exist an open cone $\calC$ and a constant $\rho>0$ such that for every 
 $x_0\in \partial\Omega$ there exist an open neighborhood $V$ and a rotation $R$ with the property that  for all  $x\in V\backslash\Omega$ we have $\Omega \cap (x +R\,\calC_\rho)=\emptyset$. 
\end{itemize}
\end{definition}
According to Proposition 9.1.4 in \cite{Agranovich-2015} (also referenced in \cite{Grisvard}), a bounded domain $\Omega\subset\R^3$ has the uniform cone property if and only if $\Omega$ has a Lipschitz continuous boundary. Furthermore, the uniform interior cone condition is equivalent to the uniform exterior cone condition (for the definition of the uniform interior cone condition, see \cite{Grisvard}).

For a function $u=(u^1,u^2,u^3)^T:\Omega\to\R^3$, the  differential  $\mathrm{D} u$ is given by
\begin{align*}
 \mathrm{D} u
 =\begin{pmatrix}
 \mathrm{D} u^1\\
 \mathrm{D} u^2\\
  \mathrm{D} u^3
 \end{pmatrix}
 \in \R^{3\times 3}\,,
\end{align*}
with $(\mathrm{D} u^k)_{l}= \partial_{x_l} u^k$ for $1\leq k\leq 3$ and $1\leq \ell\leq 3$ and with $\mathrm{D} u^k\in \R^{1\times 3}$, while for scalar functions $u:\Omega\to\R$ we also define $\grad u:=(\mathrm{D} u)^T\in \R^{3\times 1}$. For a vector field $w:\Omega\to\R^3$, the divergence and the  curl are  given as
\begin{align*}
\dyw w=\sum_{i=1}^ 3 w^{i}_{, x_i},\qquad 
\curl w =\big(w^{3}_{,x_2}-w^{2}_{,x_3},w^{1}_{,x_3}-w^{3}_{,x_1},w^{1}_{,x_2}-w^{2}_{,x_1}\big)\,.
\end{align*}
For tensor fields $Q:\Omega\to\R^{3\times 3}$, $\Curl Q$ and $\Dyw Q$ are defined row-wise: 
\begin{align*}
 \Curl Q=\begin{pmatrix} \curl Q^1\\
             \curl Q^2 \\
             \curl Q^3    
            \end{pmatrix}
\in  \R^{3\times 3}\,,\qquad 
\text{ and } 
 \Dyw Q=\begin{pmatrix}
  \dyw Q^1
  \\ \dyw Q^2 \\
   \dyw Q^3
   \end{pmatrix}\in \R^3\,,
\end{align*}
where $Q^i$ denotes the $i$-th row of $Q$. With these definitions, for $u:\Omega\to \R^3$ we have consistently $ \Curl\mathrm{D} u=0\in \R^{3 \times 3}$ and $\Dyw \mathrm{D} u=\Delta u\in \R^3$.  

The Sobolev spaces \cite{adamssobolev,Giraultbook} used in this paper are
\begin{align*}
& H^1(\Omega)=\{u\in L^2(\Omega)\, |\,\, \mathrm{D} u\in L^2(\Omega)\}\,, \qquad \|u\|^2_{H^1(\Omega)}:=\|u\|^2_{L^2(\Omega)}+\|\mathrm{D} u\|^2_{L^2(\Omega)}\, ,\\[2ex]
&H({\curl};\Omega)=\{v\in L^2(\Omega;\R^d)\, |\,\, \curl v\in L^2(\Omega)\}\,,\qquad\|v\|^2_{H({\rm curl};\Omega)}:=\|v\|^2_{L^2(\Omega)}+\|{\rm curl}\, v\|^2_{L^2(\Omega)}\, ,\\[2ex]
&H(\dyw;\Omega)=\{v\in L^2(\Omega;\R^d)\, |\,\, \dyw v\in L^2(\Omega)\},\qquad
\|v\|^2_{H({\dyw};\Omega)}:=\|v\|^2_{L^2(\Omega)}+\|\dyw v\|^2_{L^2(\Omega)}\,.
\end{align*}
Spaces for tensor valued functions are denoted by $H(\Curl;\Omega)$ and $H(\Dyw;\Omega)$. 
Moreover, $H_0^1(\Omega)$ is the completion of $C_0^{\infty}(\Omega)$ with respect to the $H^1$-norm and $H_0({\rm curl};\Omega)$ and $H_0(\dyw;\Omega)$ are the completions of $C_0^{\infty}(\Omega)$ with respect to the $H(\curl)$-norm and the $H(\dyw)$-norm, respectively. 
By $H^{-1}(\Omega)$ we denote the dual of $H_0^1(\Omega)$. Finally we define
\begin{equation} \label{space}
\begin{split}
 H(\dyw,0;\Omega)&=\{u\in H(\dyw;\Omega)\, |\,\, \dyw u=0\}\,,\\
 H(\curl,0;\Omega)&=\{u\in H(\curl;\Omega)\, |\,\, \curl u=0\}
\end{split}
\end{equation}
and set
\begin{equation} \label{space1}
\begin{split}
H_0(\dyw,0;\Omega)&=H_0(\dyw;\Omega)\cap H(\dyw,0;\Omega)\,,\\ 
 H_0(\curl,0;\Omega)&=H_0(\curl;\Omega)\cap H(\curl,0;\Omega)\,.
\end{split}
\end{equation}
%
We will obtain higher regularity of the solution for the system \eqref{Main_nonlinear} in the Besov (Nikol'skii-Besov) space $B^{m+\sigma}_{2,\infty}(\Omega)$ of differentiable functions up to order $m+\sigma$, where $m\in\mathbb{N}_0$ and $\sigma\in (0,1)$ (for a precise definition we refer to \cite{Tri78-1} and \cite{Nikolski75}). When $\Omega$ has a Lipschitz boundary, then the following identity holds true (\cite{Nikolski75,Tri78-1}) 
\begin{equation} \label{B_space}
\begin{split}
 B^{m+\sigma}_{2,\infty}(\Omega)&=\{u\in L^2(\Omega)\, |\,\, \|u\|_{\mathcal{N}^{m+\sigma}(\Omega)} < +\infty\}\,,\\[1ex]
\|u\|^2_{\mathcal{N}^{m+\sigma}(\Omega)}&= \|u\|^2_{L^2(\Omega)}+\sum\limits_{|\alpha|=m}\,\,\sup_{\substack{\eta>0\\ h\in \R^3\\0<|h|<\eta}}\int_{\Omega_{\eta}}|h|^{-2\sigma}|D^{\alpha}u(x+h)-D^{\alpha}u(x)|^2\,\di x\,,
\end{split}
\end{equation}
where $\Omega_{\eta}=\{x\in\Omega\, |\,\, \mathrm{dist}(x,\partial\Omega)>\eta\}$. It is worth emphasizing that, thanks to the family of inner variations proposed in this article (see Subsection 4.1), in the proof of the main theorem (Theorem 4.1), we will be able to estimate the second component of the norm \eqref{B_space} for the displacement $u$ ($\sigma=\frac{1}{2}$, $m=1$) and for the microdistortion $P$ ($\sigma=\frac{1}{2}$, $m=0$).

Additionally, for arbitrary $\epsilon>0$ and $s\in(0,\infty)\backslash\N$, the following embeddings are continuous: 
\begin{equation}
\label{conemb}
B^{s+\epsilon}_{2,\infty}(\Omega)\subset H^s(\Omega)\subset B^{s}_{2,\infty}(\Omega)\,,\end{equation}
where $H^s(\Omega)$ denotes the Sobolev-Slobodeckij space of order $s$.\\[1ex]
Based on the results from Section $3.3$ in \cite{Giraultbook} (Corollary $3.4$), see also  \cite[Theorem 5.3]{BaPaSc16}, the following version of the Helmholtz decomposition will be used:

\begin{theorem}[Helmholtz decomposition]
\label{thm:helmholtz-dec}
Let $\Omega\subset\R^3$ be a bounded domain with a Lipschitz boundary. Then
\[
 L^2(\Omega;\R^3)= \mathrm{D} H_0^1(\Omega)\oplus H(\dyw,0;\Omega)
\]
and hence, for every $p\in L^2(\Omega;\R^3)$ there exist a unique $v\in H_0^1(\Omega)$  and 
$q\in H(\dyw,0;\Omega)$ such that 
$p=\mathrm{D} v + q$. 
\end{theorem}
The Helmholtz decomposition theorem has as an immediate consequence (see for instance \cite{KON2023global} for a standard proof for Lipschitz domains) \begin{proposition} 
 \label{prop:helmholtz-curl}
 Let $\Omega\subset\R^3$ be a bounded domain with a Lipschitz boundary and 
let $p\in H_0(\Curl;\Omega)$ with $p=\mathrm{D} v + q$, where $v\in H_0^1(\Omega)$ and $q\in H(\Dyw,0;\Omega)$ are given according to the Helmholtz decomposition. 
 Then $\mathrm{D} v\in H_0(\Curl,0;\Omega)$ 
 and $q\in H(\Dyw,0;\Omega)\cap H_0(\Curl;\Omega)$.
\end{proposition}

The next embedding theorem is for instance proved in  \cite[Sections 3.4, 3.5]{Giraultbook} for $C^{1,1}$-domains and in \cite[Lemma 3.2]{Pauly2019} for convex domains:
\begin{theorem}[Embedding Theorem]
\label{thm:embedding}
 Let $\Omega\subset \R^3$ be a bounded convex domain or a bounded domain with a  $C^{1,1}$-smooth boundary 
 $\partial\Omega$. Then 
 \begin{align*}
  H(\curl;\Omega)\cap H_0(\dyw;\Omega)\subset H^1(\Omega),\qquad\qquad 
  H_0(\curl;\Omega)\cap H(\dyw;\Omega)\subset H^1(\Omega)
 \end{align*}
 and there exists a constant $C>0$ such that for every $p\in  H(\curl;\Omega)\cap H_0(\dyw;\Omega)$ or\\ $p\in H_0(\curl;\Omega)\cap H(\dyw;\Omega)$ we have
 \begin{equation*}
   \norm{p}_{H^1(\Omega)}\leq C(\norm{p}_{H(\curl;\Omega)} + \norm{p}_{H(\dyw;\Omega)})\,.  
 \end{equation*}
\end{theorem}
\noindent
A version of this result for Lipschitz domains showing  $H^{\frac{1}{2}}(\Omega)$-regularity can be found in \cite{Costabel90}. 
\section{Assumptions on the energy density $W$}
We consider next suitable structure assumptions on the Lagrangian energy density $W$. Let $W:\overline\Omega\times \R^{3\times 3}\times  \R^{3\times 3}\times  \R^{3\times 3}\to [0,\infty)$ be a Carath{\'e}odory-function that satisfies
\begin{itemize}
\item[(W1)] $W$ is Lipschitz in $\overline\Omega$ and of quadratic growth in the other variables. To be more precise,  there exist constants $L\geq 0$, $c_1\geq 0$ such that for all $x,y\in \overline\Omega$ and all $\bbQ=(\bbQ_1,\bbQ_1,\bbQ_3)\in  \R^{3\times 3}\times  \R^{3\times 3}\times  \R^{3\times 3}$ we have
\begin{align}
\label{W1}
\abs{W(x,\bbQ) - W(y,\bbQ)}&\leq L\abs{x-y}(1+ \norm{\bbQ}^2),\\
\label{W1_1}
W(x,\bbQ)&\leq c_1(1 + \norm{\bbQ}^2).
\end{align}
Moreover, for all $x\in \overline\Omega$, $W(x,0)=0$. 
\item[(W2)] For every $x\in\overline\Omega$ we have $W(x,\cdot)\in C^1( \R^{3\times 3}\times  \R^{3\times 3}\times  \R^{3\times 3})$ and there is a constant $c_2\geq 0$ such that for all $x\in  \overline\Omega$ and $\bbQ=(\bbQ_1,\bbQ_1,\bbQ_3)\in \R^{3\times 3}\times  \R^{3\times 3}\times  \R^{3\times 3}$  the following condition is fulfilled
\begin{align}
\label{der_W}
\norm{\rmD W(x,\bbQ)}\leq c_2(1 + \norm{\bbQ})\,.
\end{align}
\item[(W3)] The functional  $\calW:H_0^1(\Omega;\R^3) \times H_0(\Curl;\Omega)\to\R$ defined by
\begin{equation}
\label{functionalW}
 \calW(u,P)=\int_\Omega W(x,\mathrm{D} u,P,\Curl P)\,\di x   
\end{equation} is uniformly convex, i.e.\ there exists a constant $\kappa>0$ such that for all $u,v\in H_0^1(\Omega;\R^3)$ and all $P,Q\in H_0(\Curl;\Omega)$ we have
\begin{align}
\label{uconvex}
\calW(v,Q) - \calW(u,P)\,\geq &\,\,\Big\langle \rmD \calW(u,P),
\left(\begin{smallmatrix} v-u\\Q-P\end{smallmatrix}\right)\Big\rangle\nn\\[1ex] 
&+ 
\kappa \big(\norm{v-u}_{H^1(\Omega)}^2 + \norm{Q-P}^2_{H(\Curl;\Omega)}\big)\,.
\end{align}
\end{itemize}
Conditions (W1) and (W2) and the mean value theorem imply that there exists a constant $c>0$ such that with $L$ from (W1) and for all $x,y\in \overline\Omega$ and $\bbQ,\wt\bbQ\in  \R^{3\times 3}\times  \R^{3\times 3}\times  \R^{3\times 3}$ it holds
\begin{align}
\label{usefulest:W}
\abs{W(x,\bbQ) - W(y,\wt\bbQ)} \,\leq\, L\abs{x-y}(1+\norm{\bbQ}^2) + c\Big(1 + \norm{\bbQ} + \norm{\wt\bbQ}\Big)\norm{\bbQ - \wt\bbQ}\,.
\end{align}
\begin{remark}(A nonlinear example)\\
Set $\bbQ=(F,P,C)\in\R^{3\times 3}\times \R^{3\times 3}\times \R^{3\times 3}$, where $F$ and $C$ are place holders for $\mathrm{D} u$ and  $\Curl P,$ respectively, and let us consider the free energy function in the following form (skipping the material tensor for simplicity and without loss of generality)
\begin{align}
\label{nonl_freeener}
W(\bbQ)=\,\frac{1}{2}\norm{\sym(F-P)}^2+\frac{1}{2}\norm{\sym P}^2+\frac{1}{q}\norm{C}^q+\frac{1}{2}\norm{C}^2\,,
\end{align}
where $1<q<2$. Then using Young's inequality we obtain
\begin{align}
\label{nonl_freeener_1}
W(\bbQ)\leq \norm{F}^2+\frac{3}{2}\norm{P}^2+\frac{2-q}{2}\Big(\frac{1}{q}\Big)^{\frac{2}{2-q}}+\frac{q}{2}\norm{C}^2+\frac{1}{2}\norm{C}^2
\end{align}
and condition \eqref{W1_1} is satisfied. In addition, 
\begin{align}
\label{gredient_W}
\norm{\rmD W(\bbQ)}\leq \norm{\sym(F-P)}+\norm{\sym P}+\norm{C}^{q-1}+\norm{C}
\end{align}
and again by Young's inequality,  condition (W2) is satisfied as well. 
Now, let us analyze condition (W3) for the energy \eqref{nonl_freeener}. This involves considering the following expression for all $\displaystyle{u,v\in H_0^1(\Omega;\mathbb{R}^3)}$ and all $P,Q\in H_0(\text{Curl};\Omega)$
 \begin{align}
\label{nonl_uconvex}
\calW(v,Q) &- \calW(u,P)- \Big\langle \rmD\calW(u,P),
\left(\begin{smallmatrix} v-u\\Q-P\end{smallmatrix}\right)\Big\rangle= \nn\\[1ex] 
&\hspace{2ex}\frac{1}{2}\int_{\Omega}\Big(\big\langle\sym(\mathrm{D}(v-u)-(Q-P)),\sym(\mathrm{D}(v-u)-(Q-P))\big\rangle\nn\\[1ex]
&\hspace{2ex}+\big\langle\sym (Q-P),\sym (Q-P)\big\rangle
+ \big\langle\Curl (Q-P),\Curl (Q-P)\big\rangle\Big)\,\di x\\[1ex]
&\hspace{2ex}+ \int_{\Omega}\Big(\frac{1}{q}\norm{\Curl Q}^q-\frac{1}{q}\norm{\Curl P}^q
- \big\langle\norm{\Curl P}^{q-2}\Curl P,\Curl (Q-P)\big\rangle\Big)\,\di x\,.\nn
\end{align}
The real function $f(x)=x^q$ for $x\geq 0$ and $1<q<2$ is convex and differentiable, hence the last integral on the right-hand side of \eqref{nonl_uconvex} is non-negative. The incompatible Korn's inequality  \cite{OptimalKMSIneq,lewintan_neff_2021,NeffLP,neff2012poincare,Lewintan_korn_2022} implies: there is a constant $\tilde{c}>0$ such that 
\begin{equation}
\label{coercive1}
\|P\|^2_{L^2(\Omega)}\leq \tilde{c}\,\big( \|\sym  P\|^2_{L^2(\Omega)}+\|\Curl P\|^2_{L^2(\Omega)}\big)
\end{equation}
for all $P\in {\rm H}_0(\Curl;\Omega)$. Applied to the first integral on the right-hand side of \eqref{nonl_uconvex} it implies that
 \begin{align}
\calW(v,Q) - \calW(u,P)-\Big\langle \rmD\calW(u,P),
\left(\begin{smallmatrix} v-u\\Q-P\end{smallmatrix}\right)\Big\rangle\geq \kappa \big(\norm{v-u}_{H^1(\Omega)}^2 + \norm{Q-P}^2_{H(\Curl;\Omega)}\big)
\end{align}
for some constant $\kappa>0$, hence the condition (W3) is satisfied.
\end{remark}


\begin{remark}(The relaxed micromorphic model - linear case)\\
Let $\bbQ=(\mathrm{D} u,P,\Curl P)$ and let us consider the free energy function associated with the system \eqref{Main}
\begin{align}
\label{freeener}
W(x,\bbQ)=&\,\frac{1}{2}\langle \D(x)\sym(\mathrm{D} u-P),\sym(\mathrm{D} u-P)\big\rangle
+\frac{1}{2}\big\langle\Ha(x)\sym P,\sym P\big\rangle\\[1ex]
&\hspace{2ex}+ \frac{1}{2}\big\langle\Lc(x)\Curl P,\Curl P\big\rangle\,.\nn
\end{align}
We assume that the coefficient functions $\D,\Ha$ and
 $\Lc$ in \eqref{freeener} are fourth order elasticity tensors from $ C^{0,1}(\overline\Omega;\Lin(\R^{3\times 3};\R^{3\times 3}))$  and are symmetric and positive definite in the following sense
\begin{enumerate}[(i)]
 \item For every $\sigma,\tau\in \Sym(3)$, $\eta_1,\eta_2\in \R^{3\times 3}$ and all $x\in \overline\Omega$:
 \begin{gather}
  \langle \mathbb{C}_e(x)\sigma,\tau\rangle = \langle \sigma,\mathbb{C}_e(x)\tau\rangle\,,
  \quad\qquad
  \langle \mathbb{C}_\text{micro}(x)\sigma,\tau\rangle = 
  \langle \sigma,\mathbb{C}_\text{micro}(x)\tau\rangle\,,
  \nonumber\\
  \langle \mathbb{L}_c(x)\eta_1,\eta_2\rangle 
  =  \langle\eta_1, \mathbb{L}_c(x)\eta_2\rangle\,.
  \label{symmetries1}
 \end{gather}
\item There exists positive constants $C_{\mathrm{e}}$, $C_{\mathrm{micro}}$ and $L_{\mathrm{c}}$ such that for all $x\in\overline\Omega$, $\sigma\in\Sym(3)$ and $\eta\in\R^{3\times 3}$:
\begin{align}
    \label{positivedef}
  \langle\D(x)\sigma,\sigma\rangle\geq C_{\mathrm{e}}|\sigma|^2 &\,,\quad\,\, \langle\Ha(x)\sigma,\sigma\rangle\geq C_{\mathrm{micro}}|\sigma|^2\,,\quad\,\,   \langle\Lc(x)\eta,\eta\rangle\geq L_{\mathrm{c}}|\eta|^2\,. 
\end{align}
\end{enumerate}
Then we obtain that the energy \eqref{freeener} satisfies condition $(W1)$. We can also easily calculate that
\begin{align}
\label{freeener1}
\rmD W(x,\bbQ)=\begin{pmatrix}
\D(x)\sym(\mathrm{D} u-P)\\
-\D(x)\sym(\mathrm{D} u-P) + \Ha(x)\sym P\\
\Lc(x)\Curl P
\end{pmatrix}
\end{align}
and the condition $(W2)$ holds. Using the form \eqref{freeener1} of $\rmD W$ and again the incompatible Korn's inequality \eqref{coercive1}, we obtain that the condition $(W3)$ is also met.\\
\end{remark}
\noindent
Let us now assume that $f\in H^{-1}(\Omega)$ and $M\in \big(H_0(\Curl;\Omega)\big)^{\ast}$. We define the functional $\calE: H_0^1(\Omega;\R^3)\times H_0(\Curl;\Omega)\to\R$ by the formula 
\begin{equation}
    \label{functional}
    \calE(v,Q):=\calW(v,Q)-\langle f,v\rangle - \langle M,Q\rangle\,,
\end{equation}
where $\calW$ and $W$ satisfy (W1)-(W3). Then, the direct method of the calculus of variations and the uniform convexity imply the existence of exactly one minimizer $(u,P)\in H_0^1(\Omega;\R^3)\times H_0(\Curl;\Omega)$ of $\calE$ (see for example \cite{dacorogna2007direct} or \cite{evansbook}). Moreover, the minimizer $(u,P)$ of $\calE$ satisfies the following weak Euler-Lagrange equation
\begin{align}
\label{eqn:eulerlagrange}
\int_\Omega \Big\langle\rmD W(\mathrm{D} u,P,\Curl P),\left[ \left( 
\begin{smallmatrix}
\mathrm{D} v\\Q\\ \Curl Q
\end{smallmatrix}
 \right)\right]\Big\rangle\,\di x
  =\langle f,v\rangle + \langle M,Q\rangle
\end{align}
for all $(v,Q)\in H_0^1(\Omega;\R^d)\times H_0(\Curl;\Omega)$. Moreover, by choosing $Q=0$ and $v=0$ in the uniform convexity condition \eqref{uconvex}, using equation \eqref{eqn:eulerlagrange} and the weighted Young's inequality, we are able to show that there is a constant $c>0$ such that for all $f\in H^{-1}(\Omega)$ and $M\in \big(H_0(\Curl;\Omega)\big)^{\ast}$ the corresponding minimizing pair $(u,P)$ satisfies the estimate
\begin{align}
\label{est:minimnorm}
 \norm{u}_{H^1(\Omega)} +\norm{P}_{H(\Curl;\Omega)}\leq c\big(\norm{f}_{H^{-1}(\Omega)} + \norm{M}_{\big(H_0(\Curl;\Omega)\big)^{\ast}}\big)\,.\\\nn
\end{align}
\begin{remark}
If $W$ is in the form specified by equation \eqref{freeener}, then equation \eqref{eqn:eulerlagrange} serves as a weak formulation of the relaxed micromorphic model described in equation \eqref{Main}. It is important to note that for a weak solution $(u,P)$ of the system defined in equation \eqref{Main} to exist, it is sufficient to assume that $\D,\Ha,\Lc\in L^\infty(\Omega;\Lin(\R^{3\times 3};\R^{3\times 3}))$ satisfy the symmetry and positivity properties outlined in equations \eqref{symmetries1} and \eqref{positivedef} (for more detailed information, please refer to \cite{KON2023global}). 
\end{remark}
\section{Global regularity on Lipschitz domains}
\label{sec:globreg}
The main goal of this section is to prove the following regularity theorem
\begin{theorem}
\label{thm:globlipnonlin}
 Let $\Omega\subset \R^3$ be a bounded domain with a Lipschitz boundary (in the sense of graphs). Suppose that $\calW$ and $\calE$ defined in \eqref{functionalW} and \eqref{functional} fulfil assumptions (W1)-(W3). Moreover, assume that $f\in L^2(\Omega;\R^3)$  and $M\in H(\Dyw;\Omega)$. Then for every minimizer $(u,P)\in H_0^1(\Omega)\times H_0(\Curl;\Omega)$ of $\calE$ we have 
 \begin{align}
  u\in B^{3/2}_{2,\infty}(\Omega),\,\,\quad P\in B^{1/2}_{2,\infty}(\Omega),\,\,\quad \Curl P\in B^{1/2}_{2,\infty}(\Omega)
 \end{align}
 and there exists a constant $C>0$ (independent of $f$ and $M$) such that 
 \begin{align}
  \norm{u}_{B^{3/2}_{2,\infty}(\Omega)} + \norm{P}_{B^{1/2}_{2,\infty}(\Omega)} + \norm{\Curl P}_{B^{1/2}_{2,\infty}(\Omega)} \leq C\,(\norm{f}_{L^2(\Omega)} + \norm{M}_{H(\Dyw;\Omega)})\,.
 \end{align}
\end{theorem}

The proof is based on an argument using a difference quotient that relies on inner variations. This argument will be presented in Section \ref{sec:proofthmgloblipnonlin}. In the following section, we introduce appropriate families of inner variations and examine their fundamental properties. The proposed method is especially suited for the nonlinear case treated here.

\subsection{Inner variations and preliminary results}

Let $\Omega\subset\R^3$ be a bounded domain that satisfies the uniform exterior cone condition with respect to an open cone $\calC_\rho$. For a point $x_0\in\partial\Omega$ we  associate a family $\bbT(x_0)$ of admissible inner variations in the following way: 

According to the uniform cone condition there exists a constant $r>0$ and a rotation $R$ such that for all $x\in \partial \Omega\cap B_r(x_0)$ we have $(x+ R\calC_\rho)\cap\Omega =\emptyset$. Without loss of generality, in the following we assume that $R=\idd$. Choose a cut-off function $\varphi\in  C_0^\infty(\R^3)$ with $0\leq\varphi\leq 1$, $\supp\varphi\subset B_r(x_0)$ and $\varphi\big|_{B_{r/2}(x_0)}=1$. For $h\in \R^3$ we define 
\begin{align}
\label{InerVar}
 T_h:\R^3\to\R^3, \quad\quad\,T_h(x):= x + \varphi(x) h\,.
\end{align}
Obviously, we have
\begin{align}
\label{IVP1}
 \mathrm{D} T_h(x)= \idd + h\otimes \grad \varphi(x),\qquad \det \mathrm{D} T_h(x)= 1 + \langle h,\grad \varphi(x)\rangle. 
\end{align}
Then for $\abs{h}\leq \frac12\norm{\grad \varphi}_{L^\infty(\R^3)}^{-1}=:\delta$ we have 
\begin{align}
\label{IVP2}
 \abs{\det\mathrm{D} T_h(x)} = \abs{1+\langle h,\grad\varphi(x)\rangle}\geq 1-\abs{h}\abs{\grad\varphi(x)}\geq \frac{1}{2}>0
\end{align}
and if $x_1+h\varphi(x_1)=x_2+h\varphi(x_2)$, then $|x_1-x_2|\leq |h|\|\grad\varphi\|_{L^{\infty}(\Omega)}|x_1-x_2|$ which implias that $x_1=x_2$ for $|h|<\delta$. The above considerations imply that for $|h|<\delta$, the mapping $T_h:\R^3\to\R^3$ is a diffeomorphism with $T_h(B_r(x_0))=B_r(x_0)$ and $T_h(x)=x$ for $x\in \R^3\backslash B_r(x_0)$. Indeed (by contradiction), let $x\in B_r(x_0)$ and $y=T_h(x)\notin B_r(x_0)$, hence $x\neq y$. Additionally, notice that $T_h(x)=T_h(y)$ and the one to one property of $T_h(\cdot)$ implies $x=y$, which leads to a contradiction.\\
Moreover, from elementary algebra we get
\begin{align}
\label{invers}
(\mathrm{D} T_h)^{-1}(x)&=\idd-\frac{1}{1+\langle h,\grad\varphi(x)\rangle}h\otimes\grad\varphi(x)\,,\\[1ex] 
\det(\mathrm{D} T_h^{-1})(x)&=1 - \frac{\langle h,\grad\varphi(x)\rangle}{1 + \langle h,\grad\varphi(x)\rangle}\,.
\end{align}
Let us choose $h_0:=\min\{\delta,\rho\}$. Then for all $h\in \calC_{h_0}$ the mappings $T_h:\R^3\to\R^3$ are diffeomorphisms with
\begin{align}
\label{Th:mappingprop}
 T_h(x)=x \text{ for }x\in \R^3\backslash B_r(x_0),\quad \quad
 T_h(\R^3\backslash\Omega)\subset  \R^3\backslash\Omega \qquad \text{and} \quad  
 T_h^{-1}(\Omega)\subset\Omega\,.
\end{align}
The second statement is a consequence of the cone property and it implies in particular that for $x\in\partial\Omega$ we have $T_h(x)\notin\partial\Omega$. This property is essential for the construction of admissible test functions, see the proofs of Lemma \ref{lem:thu} and  Lemma \ref{lem:thcurl}.\\
We call the family $\bbT(x_0):=\{T_h\, |\,\, h\in \calC_{h_0}\}$ an admissible family of inner variations (or diffeomorphisms) associated with $x_0\in \partial\Omega$. To shorten the notation  we will write $S_h:=T_h^{-1}$ for $T_h\in \bbT(x_0)$. The following uniform estimate is valid:
\begin{align}
\label{est:Th}
 \sup_{h\in \calC_{h_0}}\Big(\norm{\det\mathrm{D} T_h}_{L^\infty(\R^3)} +
 \norm{\det(\mathrm{D} T_h)^{-1}}_{L^\infty(\R^3)}& + \norm{ T_h}_{W^{1,\infty}(\R^3)} \nn\\&+ 
 \norm{ T_h^{-1}}_{W^{1,\infty}(\R^3)}
 \Big)=:C_{x_0}<\infty\,.
\end{align}
\begin{lemma}
\label{lem:thu}
 Let $\Omega \subset\R^3$ be a bounded domain with the uniform cone property and let $x_0\in \partial\Omega$. Let furthermore $\bbT(x_0)$ be the family of admissible diffeomorphisms associated with $x_0$. For $u\in H_0^1(\Omega)$ let $\wt u$ be its extension by $0$ to $\R^3$. Then for every $h\in \calC_{h_0}$ the operator 
 $\tau_h: H_0^1(\Omega)\to H_0^1(\Omega)$ defined by $\tau_h u:= (\wt u\circ T_h)\big|_\Omega$ is well defined, linear and $\sup\limits_{h\in \calC_{h_0}}\norm{\tau_h}<\infty$, where $\norm{\tau_h}$ is understood to be the norm of the operator $\tau_h$.
\end{lemma}
\begin{proof}
Let $u\in H_0^1(\Omega)$ and $\wt u$ the extension with zero, i.e.\
$$\wt u(x)=\begin{cases}
u(x) & \text{for } x\in\Omega\,, \\
0 & \text{for }x\in\Omega\setminus\R^3
\end{cases}.$$
Clearly, $\wt u\in H^1(\R^3)$ and $\norm{u}_{H^1(\Omega)}=\norm{\wt u}_{H^1(\R^3)}$. Fix $h\in \calC_{h_0}$ and define $\tau_h u:= (\wt u\circ T_h)\big|_\Omega$. Then $\tau_hu\in H^1(\Omega)$ and $\tau_h u\big|_{\partial\Omega}=0$. The latter is an immediate consequence of the mapping properties of $T_h$ stated in \eqref{Th:mappingprop}. Hence, $\tau_h:H_0^1(\Omega)\to H_0^1(\Omega)$ is well defined. The uniform bound for the operator norm of $\tau_h$ with respect to $h\in \calC_{h_0}$ now follows from transformation formulas for integrals, the chain rule and the uniform estimate \eqref{est:Th}. 
\end{proof}
Lemma \ref{lem:thu} demonstrates that for $v \in H_0^1(\Omega)$, the function $\tau_hv$ is a valid test function for the weak formulation \eqref{eqn:eulerlagrange} if $h \in \mathcal{C}_{h_0}$.  In particular it respects the vanishing Dirichlet condition on $\partial\Omega$.

Now, let us proceed to create valid test functions for the microdistortion tensor $P$. To do this, we will introduce the following notation for the extension by zero: let $\Omega\subset\R^3$ be a bounded domain. The extension operator $\calF:L^2(\Omega)\to L^2(\R^3)$, defined by
 \begin{align*}
  \calF(p)(x)=\begin{cases}
           p(x)&\text{ for }x\in \Omega\\
           0 &\text{ otherwise}
          \end{cases}\, ,
 \end{align*}
is linear and bounded. Moreover we have the following extension property for vector fields  $p:\Omega\to\R^3$ from $H_0(\curl;\Omega)$: 
\begin{lemma}
\label{lem:extension}
Let $\Omega\subset\R^3$ be a bounded domain with Lipschitz boundary. 
 The extension operator $\calF:H_0(\curl;\Omega)\to H(\curl;\R^{3})$ 
is well defined and $\norm{\calF(p)}_{H(\curl;\R^3)} = \norm{p}_{H(\curl;\Omega)}$. 
\end{lemma}
\begin{proof}
Let $p\in H_0(\curl;\Omega)$. In the distributional sense we have for all $\phi\in C_0^\infty(\R^3,\R^3)$
\begin{align}
 \langle\curl (\calF( p)),\phi\rangle&=\int_{\R^3}\langle \calF(p), \curl\phi\rangle\,\di x=\int_\Omega  \langle p,\curl\phi\rangle\,\di x=\int_\Omega \langle \curl p,\phi\rangle\,\di x \,.
\end{align}
The last identity follows by Gauss Theorem and keeping in mind that $p\in H_0(\curl;\Omega)$. This shows that $\curl(\calF(p))$ can be interpreted as a regular distribution and we have the identity $\curl(\calF(p))=\calF(\curl p)$ in $L^2(\R^3,\R^3)$. Hence, $\calF(p)\in H(\curl;\R^3)$ and the norm identity stated in the lemma holds.
\end{proof}
We next recall the definition of a Piola-type transformation $\calP_h$ for tensor fields and of a transformation ${\cal T}_h$ that interacts well with the $\Curl$-operator, 
see for example \cite{ciarlet2021mathematical,Schoberlnotes,KON23}.   
For a function $Q:\R^3\rightarrow\R^{3\times3}$ and $h\in \R^3$ we define 
\begin{align}
\label{Testfun}
\wt{\cal T}_h(Q)(x)&:=Q(T_h(x))\,\mathrm{D} T_h(x)\,.
\end{align}
Let now $\Omega \subset\R^3$ be a bounded domain with the uniform cone property  and let $x_0\in \partial\Omega$. Let furthermore $\bbT(x_0)$ be the family of admissible diffeomorphisms associated with $x_0$ and let $S_h:=T_h^{-1}$ for $T_h\in \bbT(x_0)$. Observe that thanks to \eqref{Th:mappingprop}
for  every $h\in \calC_{h_0}$ we have $S_h(\Omega)\subset\Omega$. For a function $M: \Omega\to \R^{3\times3}$ and $h\in \calC_{h_0}$ the contravariant Piola-type transformation is given by 
\begin{align}
\forall x\in \Omega\,:\qquad \calP_h(M)(x):=\det\mathrm{D} S_h(x) (M(S_h(x))(\mathrm{D} S_h(x))^{-T}\,.
\end{align}
\begin{lemma}
\label{lem:thcurl}
 Let $\Omega \subset\R^3$ be a bounded domain with the uniform cone property  and let $x_0\in \partial\Omega$. Let furthermore $\bbT(x_0)$ be the family of admissible diffeomorphisms associated with $x_0$. Then for every $h\in \calC_{h_0}$ the operators 
 \begin{align*}
  \calT_h:=\calR_\Omega\circ\wt\calT_h\circ \calF&:H_0(\Curl;\Omega)\to H_0(\Curl;\Omega),\\
  \calP_h&:H(\Dyw;\Omega)\to H(\Dyw;\Omega)
 \end{align*}
are well defined. Here, $\calR_\Omega$ is the restriction to $\Omega$ of a function defined on $\R^3$. Moreover, 
\begin{align}
\label{curl-identity2}
 \Curl \calT_h(Q)&=\det\mathrm{D} T_h \big((\Curl \calF( Q))\circ T_h\big) (\mathrm{D} T_h)^{-T}\, ,
\\
\label{div-identity2}
\Dyw \calP_h(M)&= \det\mathrm{D} S_h\, (\Dyw M)\circ S_h\,,
 \end{align}
and there exist  constants $C_1,C_2>0$ such that for all $h\in \calC_{h_0}$,  $Q\in H_0(\Curl;\Omega)$ and $M\in H(\Dyw;\Omega)$  we have
\begin{align}
\label{lem:estnorm}
 \norm{\calT_h(Q)}_{H(\Curl;\Omega)}&\leq C_1 \norm{Q}_{H(\Curl;\Omega)}\,,\\
 \norm{\calP_h(M)}_{H(\Dyw;\Omega)}&\leq C_2 \norm{M}_{H(\Dyw;\Omega)}\,.
 \label{lem:estnormpiola}
\end{align}
The operators $\calT_h$ and $\calP_h$ are adjoint with respect to the $L^2$-inner product: for every $Q,M\in L^2(\Omega;\R^{3\times 3})$ we have
\begin{align}
\label{lem:adjointPiola}
\int_\Omega \big\langle\calT_h(Q),M\big\rangle\,\di x =\int_\Omega \big\langle Q,\calP_h(M)\big\rangle\,\di x\,.
\end{align}
\end{lemma}
\begin{proof}
First we show that $\wt\calT_h:H(\Curl;\R^3)\to H(\Curl;\R^3)$ is well defined, calculate its $\Curl$  and derive  estimate \eqref{lem:estnorm}. Let $h\in \calC_{h_0}$ and $Q\in H(\Curl;\R^3)$. Clearly, $\wt\calT_h(Q)\in L^2(\R^3)$.  Moreover, for $\psi\in C_0^\infty(\R^3;\R^{3\times 3})$ we obtain for the distributional $\Curl$ of $\wt\calT_h(Q)$ after  a change of coordinates ($y=T_h(x)$ and $x=T_h^{-1}(y)\equiv S_h(y)$):
\begin{align}
\label{weakCurl}
 \big\langle \Curl \wt\calT_h(Q),\psi\big\rangle&=\int_{\R^3}  \big\langle\wt\calT_h(Q),\Curl\psi\big\rangle\,\di x\nn\\
 &=\int_{\R^3} \big\langle Q(y)\mathrm{D} T_h\big|_{S_h(y)} ,(\Curl_x\psi (x))\big|_{x=S_h(y)}\big\rangle \det\mathrm{D} S_h(y)\,\di y.  
\end{align}
Let us recall the following identity for inner variations in connection with the $\Curl$-operator (for the proof we refer to Sch\"oberl's lecture notes \cite{Schoberlnotes})
\begin{align}
\label{curl-identity}
 \curl_x\Big( q(T(x))\,\mathrm{D} T(x)\Big) = \det(\mathrm{D} T(x))\,(\curl_{y} q)\big|_{y=T(x)} \mathrm{D} T(x)^{-T}\,.
\end{align}
With \eqref{curl-identity} we find
\begin{align}
\label{curl-identity_1}
 (\Curl_x\psi (x))\big|_{x=S_h(y)}=\frac{1}{\det\mathrm{D} S_h(y)}
 \big(\Curl_y(\psi(S_h(y))\mathrm{D} S_h(y))\big)\mathrm{D} S_h(y)^T
\end{align}
and hence
\begin{align}
\label{weakCurl1}
 \langle \Curl \wt\calT_h(Q);\psi\rangle&=
 \int_{\R^3} \big\langle Q(y)\mathrm{D} T_h\big|_{S_h(y)} ,\big(\Curl_y(\psi(S_h(y))\mathrm{D} S_h(y))\big)\mathrm{D} S_h(y)^T\big\rangle\,\di y\nn\\
 &= \int_{\R^3} \big\langle Q(y),\Curl_y(\psi(S_h(y))\mathrm{D} S_h(y))\big\rangle\,\di y\nn
 \\
 &= \int_{\R^3}  \big\langle\Curl_y Q(y) , (\psi\circ S_h(y))\mathrm{D} S_h(y)\big\rangle\,\di y\nn\\
 &=\int_{\R^3} \Big\langle\Big(\det\mathrm{D} T_h(x)\Curl_y Q(y)\big|_{y=T_h(x)}(\mathrm{D} T_h(x))^{-T}\Big),\psi(x)\Big\rangle\,\di x.
\end{align}
In \eqref{weakCurl1} we used that $Q\in H(\Curl;\R^3)$, applied again the transformation $y=T_h(x)$ and used that $\langle A,B\,C\rangle=\langle A\,C^T,B\rangle$ for $A,B,C\in \R^{3\times 3}$. 
Since by assumption the first factor in the last integral of \eqref{weakCurl1} belongs to $L^2(\R^3)$, $\Curl \wt\calT_hQ$ is a regular distribution, belongs to $L^2(\R^3)$ as well and the identity \eqref{curl-identity2} is valid on $\R^3$. For arbitrary $Q\in H_0(\Curl;\Omega)$,  Lemma \ref{lem:extension} implies $\calF(Q)\in H(\Curl;\R^3)$ and hence by the above considerations, we find $\calT_h(Q)\in H(\Curl;\Omega)$ along with the identity \eqref{curl-identity2}. Estimate \eqref{lem:estnorm} is an immediate consequence of this identity and the uniform bounds for $T_h$, c.f. \eqref{est:Th}. It remains to show that  $\calT_h(Q)$ has the correct boundary conditions on $\partial\Omega$ i.e. that $\calT_h(Q)\in H_0(\Curl;\Omega)$. Let $(Q_n)_{n\in \N}\subset C_0^\infty(\Omega)$ be an approximating sequence of $Q\in H_0(\Curl;\Omega)$ with respect to the $H(\Curl;\Omega)$-norm. Then $\calT_h(Q_n)\in C^\infty(\overline\Omega)$ and for the support we find
\begin{align}
 \supp(\calT_h(Q_n))\subset  \supp \big(Q_n\circ T_h\big) 
= T_h^{-1}(\supp Q_n). 
\end{align}
Since  $T_h^{-1}(\supp Q_n)$ is a compact subset of the open set $T_h^{-1}(\Omega)$ and since  
 $T_h^{-1}(\Omega)\subset\Omega$ by \eqref{Th:mappingprop}, $ \supp(\calT_h(Q_n))$ is a compact subset of $\Omega$ and hence, $\calT_h(Q_n)\in C_0^\infty(\Omega)$. Since $\calT_h$ is linear and bounded by \eqref{lem:estnorm}, the sequence  $(\calT_h(Q_n))_{n\in \N}$ approximates $\calT_h(Q)$ with respect to the $H(\Curl;\Omega)$-norm which finally implies that  $\calT_h(Q)\in H_0(\Curl;\Omega)$. 

 With similar arguments one verifies the properties of $\calP_h$. Indeed, let $M\in H(\Dyw;\Omega)$. Clearly, $\det\mathrm{D} S_h (\Dyw M)\circ S_h\in L^2(\Omega;\R^3)$. Hence, after a transformation of coordinates and applying the Gauss Theorem we obtain for every   $\psi\in C_0^\infty(\Omega;\R^{3})$
 \begin{align*}
 \int_\Omega  \big\langle\det\mathrm{D} S_h\, (\Dyw M)\circ S_h , \psi\big\rangle\,\di y
 &= \int_{S_h(\Omega)} \big\langle\Dyw M , (\psi\circ T_h)\big\rangle \,\di x
 = -\int_{S_h(\Omega)}  \big\langle M, \mathrm{D} (\psi\circ T_h)\big\rangle \,\di x\\[1ex]
 &=- \int_{S_h(\Omega)} \big\langle M, \mathrm{D} \psi\big|_{T_h(x)}\mathrm{D} T_h\big\rangle \,\di x
 =\int_\Omega  \big\langle -\calP_h(M),\mathrm{D} \psi \big\rangle\,\di x.
 \end{align*}
 This shows that the distributional divergence of $\calP_h(M)$ is a regular distribution,  belongs to $L^2(\Omega;\R^{3\times 3})$ and that \eqref{div-identity2} is valid. Estimate \eqref{lem:estnormpiola} is an immediate consequence of this identity.
 
Let us finally verify \eqref{lem:adjointPiola}. Having in mind that $S_h(\Omega)\subset\Omega$ and hence, $\Omega\subset T_h(\Omega)$, for $Q,M\in L^2(\Omega;\R^{3\times 3})$ it follows that 
\begin{align}
\int_\Omega  \big\langle\calT_h(Q),M \big\rangle\,\di x &= \int_\Omega  \big\langle\calF(Q)\big|_{T_h(x)} \mathrm{D} T_h(x) , M(x) \big\rangle\,\di x \\
  &= \int_{T_h(\Omega)}\Big\langle \calF(Q)(y)\mathrm{D} T_h\big|_{S_h(y)}, M(S_h(y))\det\mathrm{D} S_h(y)
  \Big\rangle\,\di y\\
  &= \int_\Omega  \big\langle Q(y),\calP_h(M)(y) \big\rangle\,\di y\,,
\end{align} 
which proves the claim.
\end{proof}
\noindent
The next lemma is crucial for estimating the right hand sides of the system \eqref{eqn:eulerlagrange} in the proof of the regularity theorem.
\begin{lemma}
 \label{lem:diffquotDivCurl}
 Let $\Omega \subset\R^3$ be a bounded domain with the uniform cone property  and let $x_0\in \partial\Omega$. Let furthermore $\bbT(x_0)$ be the family of admissible diffeomorphisms associated with $x_0$. Then there exists a constant $C>0$ such that for every 
 $h\in \calC_{h_0}$, every $P\in H_0(\Curl;\Omega)$ and $M\in H(\Dyw;\Omega)$ it holds
 \begin{align}
 \label{thm:estdiffquotdivcurl}
  \abs{\int_\Omega \big\langle (\calT_h(P)-P),M \big\rangle\,\di x} = 
  \abs{\int_\Omega \big\langle P,(\calP_h(M)-M) \big\rangle\,\di x}\leq 
  C\abs{h}\norm{P}_{H(\Curl;\Omega)}\norm{M}_{H(\Dyw;\Omega)}\,.
 \end{align}
\end{lemma}
\begin{proof} 
Let $x_0\in \partial\Omega$ with corresponding cone $\calC_{h_0}$. 
Choose $R\geq 2\,h_0 +\diam (\Omega)$ so that  $\Omega+ \calC_{h_0}\Subset B_R(0)$. 
With this, for every $h\in \calC_{h_0}$ it holds $T_h(B_R(0))\subset B_R(0)$ and $T_h(x)=x$ for every $x\in \partial B_R(0)$.
For  $h\in \calC_{h_0}$ the identity $\int_\Omega \big\langle (\calT_h(P)-P),M \big\rangle\,\di x=\int_\Omega  \big\langle P,(\calP_h(M)-M) \big\rangle\,\di x$ follows from \eqref{lem:adjointPiola}. 
 
According to \cite{BEDIVAN_extension} or \cite[Section 2.1]{electrobook} there exists a linear and continuous extension operator from $ H(\Dyw;\Omega)$ to 
$ H(\Dyw;B_R(0))$. Let $\wt M\in  H(\Dyw;B_R(0))$ be this extension of $M\in H(\Dyw;\Omega)$. We obtain 
\begin{align*}
 \int_\Omega  \big\langle (\calT_h(P)-P),M \big\rangle\,\di x &= 
 \int_{B_R(0)} \big\langle (\wt\calT_h(\calF(P))-\calF(P)),\wt M \big\rangle\,\di x
 -\int_{B_R(0)\backslash\Omega}  \big\langle\wt\calT_h(\calF(P)),\wt M \big\rangle\,\di x\,.
\end{align*}
The second term on the right hand side in fact vanishes which can be seen as follows: $x\in B_R(0)\backslash\Omega $ implies that $T_h(x)\notin\Omega$ (see \eqref{Th:mappingprop}). Therefore, $\wt\calT_h(\calF(P))(x)=\calF(P)\big|_{T_h(x)}\nabla T_h(x)=0$ for $x\in B_R(0)\backslash\Omega$.

Hence, it remains to estimate the first term on the right hand side. By Lemma \ref{lem:extension} and by the Helmholtz decomposition (Theorem \ref{thm:helmholtz-dec}) there exist a unique $q\in H_0^1(B_R(0);\R^3)$ and $Q\in H(\Dyw,0;B_R(0))$ such that $\calF( P)=\mathrm{D} q + Q$. Proposition \ref{prop:helmholtz-curl} and Theorem \ref{thm:embedding}  imply that $Q\in H^1(B_R(0))$ along with the estimate 
\begin{equation}
\label{}
\norm{Q}_{H^1(B_R(0))}\leq c\,\bnorm{\calF( P)}_{H(\Curl;B_R(0))}=c\norm{P}_{H(\Curl;\Omega)}.
\end{equation}
We obtain 
\begin{align}
 \int_{B_R(0)}  \big\langle(\wt\calT_h(\calF(P))-\calF(P)),\wt M \big\rangle\,\di x
 =& \int_{B_R(0)}  \big\langle (\wt\calT_h(\mathrm{D} q)-\mathrm{D} q),\wt M \big\rangle\,\di x\nn\\[1ex]
  &+ \int_{B_R(0)} \big\langle (\wt\calT_h(Q)-Q),\wt M \big\rangle\,\di x \,.
 \end{align}
Since $Q\in H^1(B_R(0))$, the second term on the right hand side can be estimated as
\begin{align}
 \abs{\int_{B_R(0)}  \Big\langle
 (\wt\calT_h(Q)-Q),\wt M \Big\rangle\,\di x}&\leq c\abs{h}\norm{Q}_{H^1(B_R(0))}\norm{\wt M}_{L^2(B_R(0))}
\end{align}
with a constant $c>0$ that is independent of $h\in \calC_{h_0}$, $\wt M$ and $Q$.

Similar to \eqref{lem:adjointPiola}, after applying the Gauss-Theorem and 
\eqref{div-identity2}, the first term on the right hand side can be rewritten as follows
\begin{align}
 \int_{B_R(0)} \big\langle (\wt\calT_h(\mathrm{D} q)-\mathrm{D} q),\wt M \big\rangle\,\di x &= \int_{B_R(0)}  \big\langle\mathrm{D} q,(\calP_h(\wt M) -\wt M) \big\rangle\,\di x
 \nn\\
 &=\DDDS -\DDDE \int_{B_R(0)}  \big\langle q, \Dyw(\calP_h(\wt M) -\wt M) \big\rangle\,\di x
\nn\\
&=\DDDS -\DDDE \int_{B_R(0)} \Big\langle q ,\big((\det\mathrm{D} S_h)(\Dyw\wt M)\circ S_h -\Dyw\wt M\big) \Big\rangle\,\di x
\nn\\
&=\DDDS -\DDDE\int_{B_R(0)}  \big\langle (q\circ T_h - q),\Dyw\wt M \big\rangle\,\di x\,.
 \end{align}
 Since $q\in H^1(B_R(0))$ and $\Dyw\wt M\in L^2(B_R(0))$, we  arrive at
 \begin{align*}
  \abs{ \int_{B_R(0)}  \Big\langle (\wt\calT_h(\mathrm{D} q)-\mathrm{D} q),\wt M\big\rangle\,\di x }\leq c\abs{h}\norm{q}_{H^1(B_R(0))}\norm{\wt M}_{H(\Dyw;B_R(0))}\,.
 \end{align*}
 Joining the above estimates finally proves \eqref{thm:estdiffquotdivcurl}.
\end{proof}
\subsection{Proof of Theorem \ref{thm:globlipnonlin}}
\label{sec:proofthmgloblipnonlin}
Let the assumptions of Theorem \ref{thm:globlipnonlin} be satisfied and let $(u,P)\in H_0^1(\Omega;\R^3)\times H_0(\Curl;\Omega)$ be the minimizer of $\calE$. It is sufficient to discuss the regularity in the vicinity of points $x_0\in\partial\Omega$. For such points,  test functions involve finite differences that cross the boundary of $\Omega$ and one has to guarantee
that these functions are admissible (i.e.\ respect the boundary conditions). In order to obtain higher regularity in a neighborhood of points $x_0$ from the interior of $\Omega$ by  using suitable cut-off functions one can carry out the  arguments here below with test functions that have a compact support in $\Omega$ (for similar calculations see for instant \cite{KON23}). This shortens some of the arguments here below. Hence, let $x_0\in\partial\Omega$. For $T_h\in \bbT(x_0)$ we set $u_h:=\tau_h(u)$ and $P_h:=\calT_h(P)$. Then $(u_h,P_h)\in H_0^1(\Omega;\R^3)\times H_0(\Curl;\Omega)$ according to Lemma 
\ref{lem:thu}  and Lemma \ref{lem:thcurl}. The uniform convexity of $\calW$ (see (W3)) implies:
\begin{align}
\kappa (\norm{u-u_h}^2_{H^1(\Omega)} + \norm{P-P_h}^2_{H(\Curl;\Omega)} )
&\leq - \Big\langle\rmD\calW(u,P),\left(\begin{smallmatrix} 
u_h - u\\ P_h - P
\end{smallmatrix}\right)\Big\rangle + \calW(u_h,P_h) - \calW(u,P)\,.
\end{align}
Since the pair $(u,P)$ satisfies the weak Euler-Lagrange equation \eqref{eqn:eulerlagrange}, the term involving $\rmD\calW$ can be replaced so that we are left with
\begin{align}
\label{eq:regest07}
\kappa (\norm{u-u_h}^2_{H^1(\Omega)} + \norm{P-P_h}^2_{H(\Curl;\Omega)} )
\leq &\,
\calW(u_h,P_h) - \calW(u,P) - \langle f, u_h-u\rangle_{L^2(\Omega)}\nn\\[1ex] 
&- \langle M,P_h-P\rangle_{L^2(\Omega)}\,.
\end{align}
Let us next estimate the terms on the right hand side separately. The aim is to find upper bounds of the type $C\abs{h}$, where the constant $C$ is independent of $h$.

The energy terms: after a transformation of coordinates ($y=T_h(x)$, $x=S_h(y)$), where we also use \eqref{curl-identity2}, the first term can be rewritten as follows:
\begin{align}
\label{reg-est00}
\calW(u_h,P_h)=
\int_{T_h(\Omega)} W\Big(S_h(y),\mathrm{D}_y\wt u(y) \mathrm{D} T_h\big|_{S_h(y)} , \calF(P)(y)  \mathrm{D} T_h\big|_{S_h(y)}, \nn\\[1ex]
\qquad\qquad
(\det\mathrm{D} S_h(y))^{-1}\curl_y\calF(P)(y) \mathrm{D} S_h^T\Big)\det\mathrm{D} S_h(y)\,\di y\\[1ex]
\equiv\int_{T_h(\Omega)} W(S_h(y),\bbQ_h(y))\det\mathrm{D} S_h(y)\,\di y\,.\nn
\end{align}
Since $\Omega\subset T_h(\Omega)$ (see \eqref{Th:mappingprop}) and since $\bbQ_h(y)=0$ for $y\in T_h(\Omega)\backslash\Omega$, we find with (W1) that 
\begin{align}
\label{reg-est01}
\calW(u_h,P_h)&=\int_{\Omega} W(S_h(y),\bbQ_h(y))\det\mathrm{D} S_h(y)\,\di y\\[1ex]
&=\int_{\Omega} W(S_h(y),\bbQ_h(y))\,\di y
-\int_{\Omega} W(S_h(y),\bbQ_h(y))
\left(\frac{\langle h,\grad \varphi(x)\rangle}{1 +\langle h,\grad \varphi(x)\rangle}
\right)\Big|_{x=S_h(y)} \,\di y\,.\nn
\end{align}
The last term in \eqref{reg-est01} can be estimated with (W1) and \eqref{IVP2} as 
\begin{align}
\label{eq:regest09}
\abs{\int_{\Omega} W(S_h(y),\bbQ_h(y))
\frac{\langle h,\grad \varphi\rangle}{1 +\langle h,\grad \varphi\rangle}
\,\di y}&\leq 
c\abs{h}(1 + \norm{u}^2_{H^1(\Omega)} + \norm{P}^2_{H(\Curl;\Omega)})
\end{align}
with a constant $c$ that is independent of $T_h\in \bbT(x_0)$. We will next jointly estimate the first term of \eqref{reg-est01} and the second energy term in \eqref{eq:regest07}. For that purpose let  $\bbQ(x):=(\mathrm{D} u(x), P(x),\Curl P(x))$ for $x\in \Omega$.
Clearly, there exists a constant $c>0$ that is independent of $h$ such that for almost all $x\in\Omega$
\begin{align}
\norm{\bbQ_h(x) - \bbQ(x)}\leq c\abs{h}\norm{\bbQ(x)}\,,
\end{align}
where $\bbQ_h(\cdot)$ is defined in formula \eqref{reg-est00}. With \eqref{usefulest:W} we therefore find a constant $c_1>0$ such that 
\begin{align}
\label{main_ineq}
\int_{\Omega} W(S_h(y),\bbQ_h(y))\,\di y - \int_{\Omega} W(y,\bbQ(y))\,\di y
&
\leq c_1\abs{h}\big(1+ \norm{u}^2_{H^1(\Omega)} + \norm{P}^2_{H(\Curl;\Omega)}\big)
\end{align}
and this constant again is independent of $T_h\in \bbT(x_0)$. Joining \eqref{main_ineq} with \eqref{reg-est01} and \eqref{eq:regest09},  we have shown that there exists a constant $C>0$ such that for all $T_h\in\bbT(x_0)$ we have
\begin{align}
\label{main_ineq_2}
\calW(u_h,P_h)-\calW(u,P)\leq c\abs{h}
\big(1 + \norm{u}_{H^1(\Omega)}^2 + \norm{P}_{H(\Curl;\Omega)}^2\big)
\,.
\end{align}
Let us next discuss  the term $\langle f, u_h - u\rangle_{L^2(\Omega)}$. Since $f\in L^2(\Omega)$ and $u\in H_0^1(\Omega)$ we obtain 
\begin{align}
\abs{\langle f, u_h - u\rangle_{L^2(\Omega)}}&\leq 
\norm{f}_{L^2(\Omega)}\norm{u_h - u}_{L^2(\Omega)}\leq c\abs{h}\norm{f}_{L^2(\Omega)}\norm{\mathrm{D} u}_{L^2(\Omega)}\,.
\end{align}

Finally, the term $\langle M,P_h-P\rangle_{L^2(\Omega)}$ can be  estimated with the help of Lemma 
\ref{lem:diffquotDivCurl}.
Collecting the above estimates and keeping in mind \eqref{est:minimnorm}  we have shown that there is a constant $c>0$ such that for all $T_h\in \bbT(x_0)$ it holds
 \begin{align}
 \label{final}
  \kappa (\norm{u-u_h}^2_{H^1(\Omega)} + \norm{P-P_h}^2_{H(\Curl;\Omega)} )
 &\leq c\abs{h}\Big(1+\norm{f}_{L^2(\Omega)}^2 + \norm{M}_{H(\Dyw;\Omega)}^2\Big)\,. 
 \end{align}
Dividing inequality \eqref{final} by $|h|>0$ we obtain
 \begin{align}
 \label{final_1}
  \kappa (\abs{h}^{-1}\norm{u-u_h}^2_{H^1(\Omega)} + \abs{h}^{-1}\norm{P-P_h}^2_{H(\Curl;\Omega)} )
 &\leq c\Big(1+\norm{f}_{L^2(\Omega)}^2 + \norm{M}_{H(\Dyw;\Omega)}^2\Big)\,. 
 \end{align}
Thus, for each $x_0\in\partial\Omega$, we obtain a family $ \bbT(x_0)$ of admissible inner variations for which inequality \eqref{final_1} is satisfied with a constant $c>0$ that does not depend on $T_h\in \bbT(x_0)$. Since $\partial\Omega$ is compact it can be covered with finitely many balls $B_r(x_0)$ (with $x_0\in \partial\Omega$). Therefore, from \eqref{final_1}, we obtain the following estimate for the displacement $u$ 
 \begin{align}
 \label{final_2}
\sup_{\substack{h_0>0\\ h\in \R^3\\0<|h|<h_0}}\Big(\int_{\Omega_{h_0}}|h|^{-1}\abs{u_h(x)-u(x)}^2+|h|^{-1}\abs{\mathrm{D}u_h(x)-\mathrm{D}u(x)}^2\,\di x\Big)\,\leq\, C\,,
\end{align}
where the constant $C>0$ does not depend on $\abs{h}>0$. Recalling formula \eqref{B_space} we find for the displacement that  $ u\in B^{3/2}_{2,\infty}(\Omega)$ ($\sigma=\frac{1}{2}$ and $m=1$). Again using \eqref{final_1} we obtain  that the microdistortion $P$ satisfies the following estimate 
 \begin{align}
 \label{final_3}
\sup_{\substack{h_0>0\\ h\in \R^3\\0<|h|<h_0}}\Big(\int_{\Omega_{h_0}}|h|^{-1}\abs{P_h(x)-P(x)}^2+|h|^{-1}\abs{\Curl P_h(x)-\Curl P(x)}^2\,\di x\Big)\,\leq\, C\,,
\end{align}
hence $ P\in B^{1/2}_{2,\infty}(\Omega)$ and $\Curl P\in B^{1/2}_{2,\infty}(\Omega)$ on use of \eqref{B_space}. Therefore the proof of Theorem \ref{thm:globlipnonlin} is completed.
\begin{remark}
Theorem \ref{thm:globlipnonlin} implies the following theorem of regularity of solution for the relaxed micromorphic model \eqref{Main}.
\begin{theorem}
\label{thm:globlip}
 Let $\Omega\subset \R^3$ be a bounded domain with a Lipschitz boundary in the sense of graphs. 
 Moreover, in addition let $\D,\Ha,\Lc\in C^{0,1}(\overline\Omega; \Lin(\R^{3\times 3};\R^{3\times 3}))$. 
 Finally assume that $f\in L^2(\Omega;\R^3)$  and $M\in H(\Dyw;\Omega)$. Then for every weak solution $(u,P)\in H_0^1(\Omega)\times H_0(\Curl;\Omega)$ of the system \eqref{Main} we have 
 \begin{align}
  u\in B^{3/2}_{2,\infty}(\Omega),\,\,\quad P\in B^{1/2}_{2,\infty}(\Omega),\,\,\quad \Curl P\in B^{1/2}_{2,\infty}(\Omega)
 \end{align}
 and there exists a constant $C>0$ (independently of $f$ and $M$) such that 
 \begin{align}
  \norm{u}_{B^{3/2}_{2,\infty}(\Omega)} + \norm{P}_{B^{1/2}_{2,\infty}(\Omega)} + \norm{\Curl P}_{B^{1/2}_{2,\infty}(\Omega)} \leq \,C\,(\norm{f}_{L^2(\Omega)} + \norm{M}_{H(\Dyw;\Omega)})\,.
 \end{align}
\end{theorem}
Note that the continuous embedding \eqref{conemb} shows that for arbitrary $\epsilon>0$, the microdistortion tensor  $P\in H^{\frac{1}{2}-\epsilon}(\Omega)$. This means that the regularity of $P$ is here similar to the regularity of solutions of Maxwell's equations on Lipschitz domains \cite{Bonito_maxwell,Costabel90} (for additional regularity results concerning the Maxwell equation on various domains, we refer to \cite{Buffa_Costabel,Costabel_Dauge,Kar_Sini_maxwel}).\\
\end{remark}
\footnotesize{
\noindent
{\bf\large Acknowledgements}\\
Dorothee Knees and Patrizio Neff acknowledge support within the framework of the Priority Programme SPP 2256, "Variational Methods for Predicting Complex Phenomena in Engineering Structures and Materials," funded by the Deutsche Forschungsgemeinschaft (DFG, German Research Foundation). Specifically, D. Knees is supported in the project "Rate-independent systems in solid mechanics and their coupling with other dissipative systems" (Project-ID 441222077) and P. Neff is supported in the project "A variational scale-dependent transition scheme - from Cauchy elasticity to the relaxed micromorphic continuum" (Project-ID 440935806).}

\bibliographystyle{plain}
\begin{scriptsize}
\bibliographystyle{abbrv}

\end{scriptsize}

\end{document}